# Digital twin of wood working enterprise

S N Masaev[1,2,3,6], A N Minkin[1,4], E Yu Troyak[4] and A L Khrulkevich[1,5]

[1]*Siberian Federal University*,79 Svobodnyj Avenue, Krasnoyarsk 660041,Russian Federation
[2]Department of System Analysis and Operations Research, *Reshetnev Siberian State University of Science and Technolog*, Office L-409, 410, 31 Krasnoyarsky Rabochy Avenue, Krasnoyarsk 660037, Russian Federation
[3]*Control Systems LLC*, 86 Pavlova Street, Krasnoyarsk, 660122, Russian Federation
[4]*FSBEI HE Siberian Fire and Rescue Academy EMERCOM of Russia*, 1 Severnaya Street, Zheleznogorsk 662972, Russian Federation
[5]The Main Directorate of EMERCOM of Russia for Krasnoyarsk Territory, 68 Mira Avenue, Krasnoyarsk 660122, Russian Federation

[6]Corresponding e-mail: faberi@list.ru

**Abstract**. A review of scientific papers has shown that digital twins are very common for modeling the states of physical objects. It is relevant to consider the creation of a digital twin of an enterprise to obtain various assessments in management, taking into account the human factor. In the research, the impact of the human factor is assessed by the digital twin and staff competencies. The system of goals of staff competencies stands for: cognitive, affective, psychomotor. An enterprise model is created from the description of the set of events taking place on it. Each event is mapped to an element of staff competencies target system. The resulting set of staff competencies is estimated by a universal (integral) indicator. The dynamics of the universal indicator characterizes two modes of operation of the enterprise. The first operating mode is normal. The second mode of operation is based on the implementation of staff competencies. Taking into account the competencies of personnel by Bloom's taxonomy allows you to determine their interconnection: cognitive, affective, psychomotor. It correlation depends on the work performed by the employee and the influence of the external environment. It depends on the efforts of the employee as a learner.

## 1. Introduction
Human began to use digital twins relatively recently. Computing power of computers plays a decisive role in the beginning of the use of digital twins. Since 2010, the limitation in the computing power of computers has been overcome and digital twins have been widely used.

The digital twin is most often used to create a digital model of a physical object. For example, spare parts. Now, whole technological process chains appear in the form of a digital twin. This makes the process of designing physical objects and processes cheap. This reduces management costs.

Today, the use of a digital twin taking into account the human factor and the classical theory of management in production is not considered.

Many authors have dealt with the management of such objects (systems): V V Leont'ev [1], A Krasovsky [2], L V Kantorovich [3], A G Granberg, A G Aganbegyan, V F Krotov [4], O Ressler, M Golestani and other authors with related works.

Currently, it is relevant to consider the complete digital twin of the enterprise, taking into account the personnel management methodology.

Therefore, the purpose: to assess the competence of employees of the enterprise by Bloom's taxonomy.

We going to do a tasks for that:
- Describe the enterprise model;
- Link enterprise work and method system learning of personnel;
- Apply a universal indicator (integral indicator) to assess the system of personal goals;
- Analys is of results.

Enterprise model in the classical representation [1-4].

## 2. Methods and Materials

Step 1. Enterprise events are set $X$ . $t$ is time period is analyzed. Then a enterprise model is formed $S=\{T, X\}$, where $T=\{t:t=1,..., T_{max}\}$ – a lot of time points, $x(t)=[x^1(t), x^2(t),..., x^n(t)]^T \in X$ – $n$– vector of events. Events $x^i(t)$ - the value of financial expenses and income of the enterprise. The resulting system has sufficient representation $S=\{T, X\}$. You can perform standard control operations with it [1, 4].

Step 2. Comparison of competencies $v(t)$ with events $x(t)$ of the enterprise $S=\{T, X\}$. We form $v_i^j$ (compliance $x_j^i$ is $v_i^j$ set as 1-yes, 0-no) from $i$– from competence and $j$ event model. We get competencies set $V$ where $v(t) = \left[ v_1^1(t), v_2^j(t), ..., v_m^n(t) \right]^T \in V$.

Payment competencies of the model is limited by resources $C$, then $C(V) \leq C$. It restriction applies to all subsystems of the researched system.

Step 3. Calculation of the integral $V$ index through the correlation matrix $R_i(x)$

$$V_i(t) = R_i(t) = \sum_{j=1}^{n} |r_{ij}(t)|, \tag{1}$$

$$R_k(t) = \frac{1}{k-1} \overset{o}{V}_k^T(t) \overset{o}{V}_k(t) = \|r_{ij}(t)\|, \tag{2}$$

$$r_{ij}(t) = \frac{1}{k-1} \sum_{l=1}^{k} \overset{o}{v}^i(t-l) \overset{o}{v}^j(t-l), \quad i,j = 1,...,n, \tag{3}$$

$$V_k(t) = \begin{bmatrix} v^T(t-1) \\ v^T(t-2) \\ ... \\ v^T(t-k) \end{bmatrix} = \begin{bmatrix} v^1(t-k) & v^2(t-k) & \cdots & v^n(t-k) \\ v^1(t-k) & v^2(t-k) & \cdots & v^n(t-k) \\ ... & ... & ... & ... \\ v^1(t-k) & v^2(t-k) & \cdots & v^n(t-k) \end{bmatrix}, \tag{4}$$

where $t$ are the time instants, $r_{ij}(t)$ are the correlation coefficients of the variables $v^i(t)$ and $v^j(t)$ at the time instant $t$.

Step 4. The analysis of experimental data is performed graphically. The dynamics of the integral indicator is calculated for all periods of time

$$V = \sum_{t=1}^{T=max} \sum_{i=1}^{n} V_i(t). \tag{5}$$

### 2.1. Characteristics of the research objects

The company employs 650 people. These are administrative staff and workers. The company is engaged in the processing of round timber. Has three key business processes. The first one is logging. Second. Delivery along the Yenisei River to the plant. The third is the production of round wood products. A total of 1.2 million processes have been modeled. These processes can be studied in a separate work [5, 6].

To assess the competencies of the company's personnel, we use the method for assessing the achievement of goals at a woodworking enterprise. This method was named Bloom's Taxonomy in honor of its developer, American author B Bloom. Bloom's taxonomy is being implemented to track the level of competence of enterprise personnel.

Taxonomy is a system of educational goals of student [7]. The traditional management of goal system is Bloom's taxonomy [7-10].

Coincidence or not, the method appeared at the very peak of the Cold War. Not only strong and brave people were required, but smart and economical ones. Obviously this technique, like many others, appeared in view of this request.

For its time, it should be considered as a revolutionary technique. The usual pattern in teaching was "remember and tell at the right time." Pupils often repeated after the teacher in order to remember better. Bloom's taxonomy implies a scheme of "making statements from a set of evidence and arguments". It is much more complicated.

Bloom's taxonomy isn't the only way to assess the competence of student or employee learners. In parallel, other approaches to assessing competencies have developed through Universal Competencies (recommended by the Council of Europe) [11], Dublin Descriptors [12], European Qualifications Framework (passport of qualifications) [13], European Qualifications Framework for EU countries [14], National Qualifications Framework [15].

In 1960, Benjamin Bloom finally formed the foundations of the method [10]. The student's system of goals is characterized by competencies: cognitive, affective, psychomotor.

There are many objective criticisms of the application of Bloom's taxonomy. It is very difficult to confirm or refute the method itself, since there are no precise assessments of competencies, trainees. Learner Goal System (Bloom's taxonomy):

　　1. Cognitive. Knowledge –mental abilities of a person that he can use consciously [9].
　　2. Affective. Perception - covers issues related to the emotional component of learning (from the basic desire for readiness to receive information to the integration of beliefs, ideas and views. Act, adhere to, inquire, acknowledge, respond, help, challenge, cooperate, protect, adapt, demonstrate, distinguish, initiate, etc. Reacting Value orientations. Organization.
　　3. Psychomotor. Imitation - it mainly covers physical skills, which includes the coordination of brain and muscle activity. adapt, assemble, balance, build, combine, copy, design, produce, discover, discern, dissect, perform, imitate, imitate, manipulate, recognize, etc. Control. Accuracy. Joint. Naturalization.

Simulation is performed in the software package described in a separate work [5].

*2.2. Experiment result*

Given: $n$=1.2 million values, $X$=5 641 442 thousand rubles, control is set through Bloom's taxonomy ($V_{taxonomy}$). From the 7th period, three HR managers are hired in accordance with the directions of the company's activities to conduct the selected management at the enterprise. After 6 periods, managers are dismissed from their posts. The calculation algorithm is 417 minutes.

Calculation of the normal operating mode ($V_{basic\_mode}$) can be viewed in a separate work [6].

A table 1 shows the experiment result of estimating the control mode $V_i(t)$ through competence of Bloom's taxonomy.

**Table 1.** Regime $V_{(taxonomy)}$.

| $t$ | $V_{(taxonomy)}$ | $t$ | $V_{(taxonomy)}$ | $t$ | $V_{(taxonomy)}$ | $t$ | $V_{(taxonomy)}$ |
|---|---|---|---|---|---|---|---|
| 1 | 110.67 | 16 | 70.08 | 31 | 104.10 | 46 | 103.64 |
| 2 | 90.31 | 17 | 101.54 | 32 | 97.66 | 47 | 87.22 |
| 3 | 68.99 | 18 | 132.61 | 33 | 81.19 | 48 | 68.26 |
| 4 | 81.36 | 19 | 132.20 | 34 | 75.22 | 49 | 54.14 |
| 5 | 90.01 | 20 | 153.92 | 35 | 68.52 | 50 | 62.51 |
| 6 | 110.04 | 21 | 164.53 | 36 | 60.51 | 51 | 49.02 |

| 7 | 123.53 | 22 | 150.99 | 37 | 53.13 | 52 | 60.58 |
| --- | --- | --- | --- | --- | --- | --- | --- |
| 8 | 124.07 | 23 | 144.74 | 38 | 61.65 | 53 | 59.32 |
| 9 | 125.43 | 24 | 119.09 | 39 | 53.51 | 54 | 147.33 |
| 10 | 100.32 | 25 | 91.02 | 40 | 51.84 | 55 | 158.41 |
| 11 | 74.76 | 26 | 104.61 | 41 | 72.03 | 56 | 156.87 |
| 12 | 72.97 | 27 | 91.60 | 42 | 93.08 | 57 | 167.90 |
| 13 | 79.69 | 28 | 79.04 | 43 | 99.23 | Total | 5491.18 |
| 14 | 90.79 | 29 | 76.26 | 44 | 115.78 | | |
| 15 | 67.93 | 30 | 95.32 | 45 | 110.11 | | |

A figure 1 shows the experiment result of estimating the control mode $V_i(t)$ through Bloom's taxonomy.

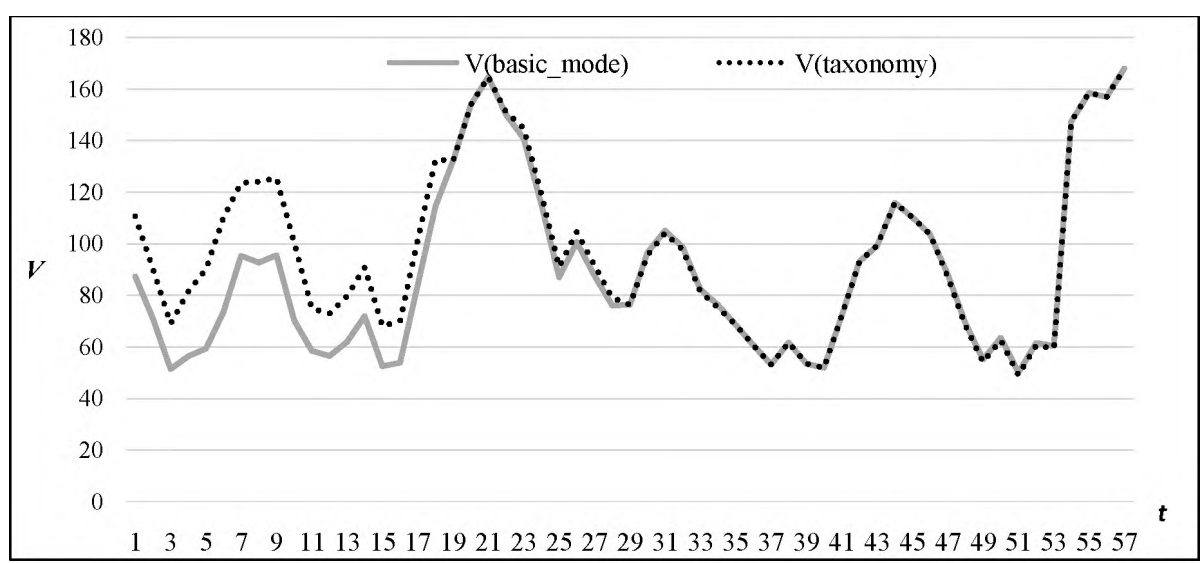

**Figure 1.** Indicator dynamics $V_i(t)$.

## 3. Results and Discussion

The installation of Bloom's taxonomy for assessing the competence of personnel in the enterprise costs 32 809 thousand rubles. Then the total costs of the enterprise for five years will amount to 5 674 251 thousand rubles.

The calculated universal indicator (integral) characterizes the possibility of introducing Bloom's taxonomy. The method requires a lot of work and resources to implement.

Assessing the competence of the enterprise personnel using Bloom's taxonomy, we obtain an assessment of the mental and physical abilities of a person. We get an idea in what period of time they were applied. Allows you to assess the need for staff training. Assess the effectiveness of retraining of workers and workers. Assess the level of training of graduates who came to work at the enterprise.

A similar approach has been used to optimally control large economic assets, such as the special economic zone of the city of Krasnoyarsk [16, 17].

## 4. Conclusion

In a simulated digital twin of a wood processing enterprise, all the main production and administrative management processes are taken into account: procurement of raw materials, transportation and logistics, loading, unloading, processing, warehouse control, sales. On a digital twin of the enterprise, HR management was assessed, taking into account its competencies by Bloom's taxonomy.

The learner's goal system (Bloom's taxonomy) has not lost its value and relevance. The obtained universal assessment (integral assessment) will allow avoiding criticism of the use of this approach. It clearly establishes the relationship between the system of goals and competencies of the student:

cognitive, affective, psychomotor. It makes it possible to apply Bloom's taxonomy to larger industrial properties.
   The research tasks were completed:
   - The model of the enterprise is characterized $S=\{T, X\}$;
   - United the work of the enterprise and the system of goals learning from Bloom's taxonomy
   $$v(t) = \left[v_1^1(t), v_2^j(t), \ldots, v_m^n(t)\right]^T \in V;$$
   - Epy universal indicator (integral indicator) was applied to assess the system of students' goals $V_{(taxonomy)}$;
   - Results analysis performed $\Delta V = V_{(taxonomy)} - V_{(basic\_mode)} = 5\,491.18 - 5\,069.93 = 421.25$.

   The purpose of the research has been achieved.